\documentclass{article}
\usepackage{amsmath, amsthm, amssymb, bm , graphicx , subfigure, array, color, float,sectsty,hyperref}
\hypersetup{
    bookmarks=true,         
    unicode=false,          
    pdftoolbar=true,        
    pdfmenubar=true,        
    pdffitwindow=false,     
    pdfstartview={FitH},    
    pdftitle={My title},    
    pdfauthor={Author},     
    pdfsubject={Subject},   
    pdfcreator={Creator},   
    pdfproducer={Producer}, 
    pdfkeywords={keywords}, 
    pdfnewwindow=true,      
    colorlinks=false,       
    linkcolor=red,          
    citecolor=green,        
    filecolor=magenta,      
    urlcolor=cyan          
}
\sectionfont{\LARGE}
\subsectionfont{\Large}
\newtheorem*{Theorem1}{Theorem 2.2 (the inequality)}
\newtheorem*{Theorem8}{Lemma 2.1}
\newtheorem*{Theorem2}{Proposition 3.1}
\newtheorem*{Theorem3}{Lemma 3.9}
\newtheorem*{Theorem4}{Proposition 3.7}
\newtheorem*{Theorem5}{Lemma 4.1}
\newtheorem*{Theorem11}{Proof of theorem \hyperlink{theorem:1.5}{1.5}}
\newtheorem*{Theorem12}{Lemma 3.4}
\newtheorem*{Theorem13}{Lemma 3.5}
\newtheorem*{Theorem14}{Lemma 3.6}
\newtheorem*{Theorem15}{Lemma 3.3}
\newtheorem*{Theorem16}{Lemma 3.2}
\newtheorem*{Theorem17}{Lemma 3.8}
\begin{document}
\setcounter{footnote}{1}
\title{A characterization of the disc through a hessian equality}
\author{Netanel Blaier\footnote{Supported by US-Israel Binational Science Foundation grant 2006136.}}
\date{\today}
\maketitle
\begin{abstract}
Let $M$ be a bounded open plane domain. Let $f$ be a continuous function on the closure of $M$, $3$-times continuously differentiable in $M$, which vanish on the boundary. Polterovich and Sodin proved that the values of $f$ cannot exceed the norm of the hessian of $f$, averaged over the entire domain $M$. In this paper we study the equality case for this inequality. We show that equality holds if and only if $M$ is a open disc and $f$ belongs to a special class of radial functions. 
\end{abstract}
\section{Introduction and summary of results} 
Let $M$ be a bounded open domain in $\mathbb{R}^2$. For every $k>0$, denote by $C^k_0(M)$ the space of all continuous functions $f:\overline{M} \rightarrow \mathbb{R}$ such that $f \in C^k(M)$ and $f$ vanishes on $\partial M$. \\\\
In their article \cite[remark 3.5]{[1]} Polterovich and Sodin proved the following inequality\footnote{In the paper \cite{[1]} the result is stated with slightly stronger assumptions on $M$ and $f$, but the proof is valid for $C^2$ functions and open domains.}: For every $f \in C^2_0(M),$ \hypertarget{eq:ineq}{}
\begin{equation}
\underset{M}{\max}|f|\leq \frac{1}{2\pi}\int_{M}||H_f||d\sigma
\end{equation}
where $H_f$ denotes the Hessian of $f$ and $|| \cdot ||$ is the operator norm.
Also of interest is the fact that this inequality was proven using the Sasaki metric and the Banach indicatrix; that is, using geometric arguments instead of the standard analytic tools. 

A question naturally arises: if equality holds in the above inequality, what can be said about the geometric nature of the domain $M$ and the function $f$ ?      

Our first result is as follows, \hypertarget{theorem:1.1}{} \\\\
\textbf{Theorem 1.1:} \textit{Let $M \subset \mathbb{R}^2$ be a bounded open domain. Let $f \in C^3_0(M)$ be a non constant function. If 
$$\underset{M}{\max}|f| = \frac{1}{2\pi}\int_{M}||H_f||d\sigma$$
then $M=D(p,R)$; that is, $M$ is a open disc.} \hypertarget{theorem:1.2}{} \\\\
\textbf{Remark 1.2:} Note that inequality (\hyperlink{eq:ineq}{1}) is clearly invariant under euclidean isometries and dilations. So it follows from the theorem above that in order to study the equality case we can restrict our investigation to the disc $M=D(0,1)$.\\\\
Let $\mathcal{S} \subseteq C^3_0(D(0,1))$ be the family of non constant solutions to the equality 
$$\underset{D(0,1)}{\max}|f| = \frac{1}{2\pi}\underset{D(0,1)}{\int} ||H_f||d\sigma$$
with the normalization $f(0,0)=1$. Denote by $\mathcal{G}$ the following set of functions
$$\mathcal{G}=\left\{g:[0,1] \rightarrow \mathbb{R}\left|\right. \: \: \mathrm{g \: \: is  \: \: C^2-smooth \: \: and}\: \: \forall t \in [0,1): -\frac{1}{t} \leq g'(t) \leq 0 \right\}$$
with normalization $\int_0^{1} e^{g(\tau)} d\tau = 1$.
For every $g \in \mathcal{G}$ we associate a function $f_g : \overline{D}(0,1) \rightarrow \mathbb{R}$ defined
by the rule $$f_g(x,y)= \int_{x^2+y^2}^{1} e^g(\tau) d\tau.$$
Our next result shows that all the normalized solutions to (\hyperlink{eq:ineq}{1}) are obtained in this manner. \hypertarget{theorem:1.3}{} \\\\
\textbf{Theorem 1.3:} \textit{(1) For every $g \in \mathcal{G}$ the function $f_g$ is a normalized solution to the equality. That is, $f_g \in \mathcal{S}$.\\
(2) For every function $f \in \mathcal{S}$ there exists a unique $g \in \mathcal{G}$ such that $f=f_g$.}\\\\
We have the following obvious properties from the classification of the $\mathcal{S}$: \\\\
\textbf{Corollary 1.4:} \textit{Let $f \in \mathcal{S}$. Then:\\
(1) $f$ is a radial function.\\
(2) $f$ is strictly monotone decreasing with respect to $r$.\\
(3) $f$ is non-negative.\\
(4) $f(x,y)=0$ if and only if $x^2+y^2=1$.\\
(5) If $|\nabla f|_{(x,y)}=0$ then $(x,y)=(0,0)$.}\\\\
Now, $\mathcal{S}$ is not empty. The most immediate example of a solution to (\hyperlink{eq:ineq}{1}) is $f(x,y)=1-(x^2+y^2)$. It turns out that this example has a very nice property: it bounds $\mathcal{S}$ from below. \hypertarget{theorem:1.5}{} \\\\
\textbf{Theorem 1.5:} \textit{Let $f \in \mathcal{S}$. Then for every $(x,y) \in \overline{D}(0,1):$ $$f(x,y) \geq 1 - (x^2+y^2)$$}
\subsection{Structure of the paper}
In article \cite{[1]} Polterovich and Sodin proved inequality (\hyperlink{eq:ineq}{1}) using the Sasaki metric. In order to understand the equality case, it is instrumental to reprove the inequality using only elementary curve properties. Section \hyperref[section2]{2} is dedicated to that purpose. In section \hyperref[section3]{3}
we study the geometry of the level sets and the behavior of $f$ and $|\nabla f|$ along the gradient flow.
This would allow us to prove theorem \hyperlink{theorem:1.1}{1.1}. In section \hyperref[section4]{4} we study the set of functions for which equality holds in (\hyperlink{eq:ineq}{1}) and prove theorems \hyperlink{theorem:1.3}{1.3} and \hyperlink{theorem:1.5}{1.5}. \\\\
\textbf{Acknowledgements.} I would like to thank my advisor Prof. Leonid Polterovich for providing me with the opportunity to work on the equality case, and also for many hours of guidance and numerous corrections and improvements.  

I would also like to thank Prof. Michael Sodin for reviewing the paper and refining some of the results,
and Prof. Fedor Nazarov for communicating to me some ideas used in this paper.
\section{Proof of the inequality}
\label{section2}
Let $f \in C^2_0(M)$. We will say that $c \in \mathbb{R}$ is a regular value (critical value) if it is a regular value (critical value) for $f|_M$. Let $c \in \mathbb{R}$ be a regular value. Then $f^{-1}(c)$ is a one-dimensional manifold without boundary. So $f^{-1}(c)$ is the disjoint union $$f^{-1}(c)=\overset{n}{\underset{i=1}{\cup}} \gamma_{i}$$ 
where every $\gamma_{i}$ is a simple closed curve. \\\\
\textbf{Notation:} For every such curve $\gamma_i$ we introduce the following notations: \\
(i) Let $L_i$ denote the length of the curve $\gamma_i$. Denote by $\gamma_{i}(s):\left[0,L_i\right] \rightarrow M$ the arclength parameterization on $\gamma_i$. \\
(ii) For every $s \in \left[0,L_i\right]$ we denote by $w_i(s)=\dot{\gamma_{i}}(s)$ the tangent vector, and by 
$\nu_i(s)$ the outward unit normal. \\
(iii) Denote by $\kappa_i(s):[0,L_i] \rightarrow \mathbb{R}$ the signed curvature of $\gamma_i$. \\\\
Then we have a simple identity relating the curvature of $\gamma_i$ to the hessian of $f$:\hypertarget{lemma:2.1}{}
\begin{Theorem8} For every $1 \leq i \leq n:$ 
$$ |\kappa_i(s)| = \frac{|<H_f w_i(s), w_i(s)>|}{|\nabla f|} .$$
\begin{proof}
Let $1 \leq i \leq n$. Denote $w = w_i$ and $\nu = \nu_i$.
The vectors $w,\nu$ form a Frenet frame on $\gamma_i$. Denote by $\nabla_w \nu$ the covariant derivative of the vector field $\nu$ in the direction of the vector field $w$. From Frenet-Serret equations in $\mathbb{R}^2$ we see that
$$|\nabla_w \nu (s)| = |\dot{\nu} (s)| = |\kappa_i(s)| |w(s)| = |\kappa_i(s)|.$$
Since $\nu$ is the unit normal $<\nu(s),\nu(s)> = 1$ for all $0 \leq s \leq L_i\:$. Therefore $$<\nabla_w \nu,\nu>=\frac{1}{2}\nabla_w<\nu,\nu> \equiv 0$$ 
and
$$\nabla_w \nu = <\nabla_w \nu , w>w + <\nabla_w \nu , \nu>\nu = <\nabla_w \nu , w>w.$$
But 
\begin{center}
$\nabla_w \nu = \nabla_w (\frac{\nabla f}{|\nabla f|})= \frac{1}{|\nabla f|} \underbrace{\nabla_w \nabla f}_{H_f} + (\mathcal{L}_w \frac{1}{| \nabla f|}) \underbrace{\nabla f}_{\lambda \cdot \nu} $ \end{center} 
which proves the statement above. 
\end{proof}
\end{Theorem8} 
$ $\\
\textbf{Notation:} Let $f \in C^2_0(M)$ and let $c \in \mathbb{R}$ be a regular value. Denote by $\beta(c, f)$ the number of connected components of $f^{-1}(c)$. Let $u \in C(\mathbb{R})$. We define the \textsl{generalized Banach indicatrix} as 
$$ B(u,f) = \overset{\infty}{\underset{-\infty}{\int}} u(c) \beta(c, f) dc.$$
\begin{Theorem1} Let $M$ be a bounded open domain in $\mathbb{R}^2$. For every $f \in C^2_0(M)$:
$$\underset{M}{\max}|f|\leq \frac{1}{2\pi}\int_{M}||H_f||d\sigma .$$
\begin{proof}
Note that if $f$ is constant on $M$ then $f \equiv 0$, and the inequality obviously holds. So we can assume that $f$ is not constant. The Banach indicatrix $B(1,f)$ bounds the left side of the inequality from above:
$$B(1,f) = \int^{\infty}_{-\infty}\beta(c,f)dc \geq \int^{\underset{M}{\sup} f}_{\underset{M}{\inf} f} 1 dc = \underset{\overline{M}}{\max} f - \underset{\overline{M}}{\min} f.$$
and since $f$ vanishes on $\partial M$ necessarily $\underset{M}{\max}|f|\leq B(1,f)$.\\\\
Let $c \in \mathbb{R}$ be a regular value. Using the same notations as in lemma \hyperlink{lemma:2.1}{2.1}, we can define the following function $$\displaystyle L(c) := \sum^{\beta(c)}_{i=1} \int_{\gamma{i}}|\kappa(s)|ds \: .$$ 
Note that 
$$L(c) = \sum^{\beta(c)}_{i=1} \int_{\gamma{i}}|\kappa(s)|ds \geq \sum^{\beta(c)}_{i=1} \int_{\gamma{i}} \kappa(s) ds = 2\pi \beta(c,f).$$
Where the last equality follows from Hopf theorem for the total curvature of a simple closed connected plane curve. Now we can integrate over all the regular values $c \in \mathbb{R}$ of $f|_M$ and get an upper bound for $B(1,f)$:
$$ B(1,f) = \frac{1}{2\pi} \int^{\infty}_{-\infty} 2 \pi \beta(c,f) dc \leq \frac{1}{2\pi} \int^{\infty}_{-\infty}  L(c) dc. $$
$L(c)$ itself is bounded from above by the operator norm of the hessian:
\begin{align*}
L(c) &= \sum^{\beta(c)}_{i=1} \int_{\gamma{i}}|\kappa(s)|ds = \sum^{\beta(c)}_{i=1} \int_{\gamma{i}} \frac{|<H_f w(s), w(s)>|}{|\nabla f|} ds \leq \sum^{\beta(c)}_{i=1} \int_{\gamma{i}} \frac{||H_f||}{|\nabla f|} ds \\ &= \int_{f^{-1}(c)} \frac{||H_f||}{|\nabla f|} ds \:.
\end{align*}
We observe that the integral on the right-hand side equals $\underset{\epsilon \rightarrow 0^+}{\lim} I_\epsilon$ where
$$\displaystyle I_\epsilon := \int_{f^{-1}(c)} \frac{||H_f||}{|\nabla f|+\epsilon} ds \: .$$  
For $\epsilon>0$ the integrand $\displaystyle \frac{||H_f||}{|\nabla f|+\epsilon}$ is defined for all $M$, not just on regular points.  Thus we can apply the smooth co-area formula to $I_\epsilon$:
\begin{align*}
\int^{\infty}_{-\infty} I_\epsilon dc &= \int^{\infty}_{-\infty} \left(\int_{f^{-1}(c)} \frac{||H_f||}{\epsilon + |\nabla f|} ds \right)dc = \int_{M} \frac{|\nabla f|}{\epsilon + |\nabla f|} ||H_f|| d \sigma \\ &\leq \int_{M} ||H_f|| d \sigma \: .
\end{align*}
Since singular values are a set of measure 0 (by Sard lemma), passing to the limit as $\epsilon \rightarrow 0^+$ we see that
\begin{align*}
\int^{\infty}_{-\infty} \left( \int_{f^{-1}(c)} \frac{||H_f||}{|\nabla f|} ds \right) dc = \int^{\infty}_{-\infty} I_0 dc \leq \int_{M} ||H_f|| d\sigma \: 
\end{align*}
\hypertarget{eq:ineq2}{}
where the integrand on the left is understood to be zero for singular values of $f$. Combining all the results above gives us a chain of estimates 
\begin{align}
\underset{M}{\max}|f| &\leq B(1,f) \leq \frac{1}{2 \pi} \int^{\infty}_{-\infty} L(c) dc \leq 
\frac{1}{2 \pi} \int^{\infty}_{-\infty} \left(\int_{f^{-1}(c)} \frac{||H_f||}{|\nabla f|} ds\right) dc \nonumber \\
& \leq \frac{1}{2 \pi}\int_{M}||H_f|| d\sigma \: .
\end{align}
\end{proof}
\end{Theorem1}
\section{The equality case : proof that M is an open disc}
\subsection{Analysis of the equality condition}
In this section we analyze the equality condition. We will show that
$$\underset{M}{\max}|f| = \frac{1}{2\pi}\int_{M}||H_f||d\sigma$$
if and only if $f$ has certain properties. Some of these properties concern the local behavior of $f$ near regular level sets, while other place a bound on the rate of descent of $f$ along the normal flow (the flow of the vector field $\pm \frac{\nabla f}{|\nabla f|}$).

We need to make the following notations: denote by $\nu = \pm \frac{\nabla f}{|\nabla f|}$ the field of outward unit normals on $M$. Let $L_\nu f$ denote the Lie derivative of $f$ along $\nu$, and let $\mu$ be the standard Lebesgue measure on $\mathbb{R}^2$. Then we have the following proposition: \\ \hypertarget{theorem:3.1}{}
\begin{Theorem2} Let $0 \not\equiv f \in C^3_0(M)$. $f$ attains equality in (\hyperlink{eq:ineq}{1}) if and only if for every regular value $c \in \mathbb{R}$:\\
(1) $f^{-1}(c)$ is a convex plane curve. \\
(2) $f$ is either non-negative or non-positive.\\
(3) $\nabla_{\nu}\nu \equiv 0$ on the level set $f^{-1}(c)$.\\
(4) $|\nabla f|$ and $L_\nu |\nabla f|$ are constant on $f^{-1}(c)$.\\
(5) $|L_\nu |\nabla f|| \leq \frac{2\pi}{Length(f^{-1}(c))}|\nabla f|$. \\
(6) $\mu(\left\{ x \in M \left| \right. {|\nabla f|}_x=0 \right\})=0$. 
\begin{proof}
On one direction, we will assume that the equality $$\underset{M}{\max}|f| = \frac{1}{2\pi}\int_{M}||H_f||d\sigma$$ 
holds and show that properties (1) to (6) follow. \\ \hypertarget{eq:eq3}{} \\
First, we note that because $\underset{M}{\max}|f| = \frac{1}{2\pi}\int_{M}||H_f|| d\sigma$ the chain (\hyperlink{eq:ineq2}{2}) must be a chain of equalities. That is, 
\begin{align}
\underset{M}{\max}|f| &= B(1,f) = \frac{1}{2 \pi}  \int^{\infty}_{-\infty} L(c) dc = 
\frac{1}{2 \pi} \int^{\infty}_{-\infty} \left(\int_{f^{-1}(c)} \frac{||H_f||}{|\nabla f|} ds\right) dc \nonumber \\
&= \frac{1}{2 \pi}\int_{M}||H_f|| d\sigma \:.
\end{align} 
Examining the leftmost equality
$$\underset{M}{\max} |f| = B(1,f) = \int^{\infty}_{-\infty}\beta(c,f)dc \geq \int^{\underset{M}{\sup} f}_{\underset{M}{\inf} f} 1 dc = \underset{\overline{M}}{\max} f - \underset{\overline{M}}{\min} f \geq \underset{M}{\max} |f| $$ 
we see that necessarily $\beta(c,f) \equiv 1$ for almost all regular values $c \in \mathbb{R}$ of $f|_M$ and $f$ is either non-negative or non-positive. This proves property (2). Note that $\beta(c,f) = 1$ for almost every regular value implies that $\beta(c,f) = 1$ for every regular value.\\\\
Looking at the next equality in (\hyperlink{eq:eq3}{3})
\begin{align}
\overset{\infty}{\underset{-\infty}{\int}} \beta(c,f) dc &= B(1,f) =  \frac{1}{2\pi} \overset{\infty}{\underset{-\infty}{\int}} L(c) dc \nonumber 
\end{align}
so for almost every regular $c \in \mathbb{R}$
\begin{align}
\frac{1}{2\pi}\int_{f^{-1}(c)}|\kappa(s)|ds = L(c) = \beta(c,f) = 1. \nonumber
\end{align}
Note that $c \in \mathbb{R}$ is a regular value and $\beta(c,f)=1$ so $f^{-1}(c)$ is a simple, closed connected  plane curve. We know that the total curvature of the curve $f^{-1}(c)$ is $2\pi$, so by Fenchel-Borsuk theorem $f^{-1}(c)$ is a convex plane curve. This proves property (1). \\\\
Also from (\hyperlink{eq:eq3}{3})
\begin{align}
\frac{1}{2\pi} \int_{f^{-1}(c)} \frac{|<H_f \dot{\gamma}(s), \dot{\gamma}(s)>|}{|\nabla f|}ds = L(c) = \frac{1}{2\pi} \int_{f^{-1}(c)}\frac{||H_f||}{|\nabla f|}ds \nonumber
\end{align}
so $||H_f|| = |<H_f \dot{\gamma}(s) ,\dot{\gamma}(s)>|$ for almost every $s \in f^{-1}(c)$. Again because of a continuity argument and Sard lemma we can drop these restrictions. So $||H_f|| = |<H_f \dot{\gamma}(s) ,\dot{\gamma}(s)>|$ for every regular value $c \in \mathbb{R}$ and every $s \in f^{-1}(c)$. From basic functional analysis that is equivalent to the following statement: for every regular value $c \in \mathbb{R}$ and for every $s \in f^{-1}(c)$, $\dot{\gamma}(s)$ is an eigenvector of $H_f$ corresponding to the eigenvalue with the largest absolute value. \\\\
Observe that if $\mathcal{R}$ is the set of regular values of $f|_M$ and $U=f^{-1}(\mathcal{R})$ then $U$ is open from Sard lemma. Denote $w = \dot{\gamma}(t)$ and $\nu = \pm \frac{\nabla f}{|\nabla f|}$ the outward unit normal. The vector fields $w,\nu$ form an orthonormal frame on $U$ i.e. $<\nu,\nu> = <w,w> \equiv 1$ and $<\nu,w>\equiv 0$ on $U$. Taking the Lie derivative
\begin{align}
0 &= \mathcal{L}_\nu <w,w> = <\nabla_\nu w , w> + <w , \nabla_\nu w> = 2<\nabla_\nu w , w> \nonumber \\
0 &= \mathcal{L}_w <\nu,\nu> = <\nabla_w \nu , \nu> + <\nu , \nabla_w \nu> = 2<\nabla_w \nu , \nu> \nonumber \\
0 &= \mathcal{L}_\nu <\nu,\nu> = <\nabla_\nu \nu , \nu> + <\nu , \nabla_\nu \nu> = 2<\nabla_\nu \nu , \nu> \nonumber
\end{align}
so $\nabla_\nu w \perp w \: , \: \nabla_w \nu \perp \nu \: ,\: \nabla_\nu \nu \perp \nu$. The vectors $w$ and $\nu$ are also eigenvectors of the Hessian, so they form an eigenbasis.
But from the definition of the Hessian
$$H_f \nu = \nabla_{\nu}\nabla f = \nabla_{\nu}(|\nabla f| \nu) = |\nabla f| \nabla_{\nu}\nu + \mathcal{L}_\nu |\nabla f| \nu$$ 
$$H_f w = \nabla_w\nabla f = \nabla_w(|\nabla f| \nu) = |\nabla f| \nabla_w \nu + \mathcal{L}_w |\nabla f| \nu$$
so $\mathcal{L}_w |\nabla f| \equiv 0$ and $\nabla_{\nu}\nu \equiv 0$. Taking the Lie derivative from the inner product
\begin{eqnarray*}
0 = \mathcal{L}_\nu <w,\nu> = <\nabla_\nu w , \nu> + <w , \nabla_\nu \nu> = <\nabla_\nu w , \nu>.
\end{eqnarray*}
This shows that $\nabla_\nu w \perp \nu$ also, so $\nabla_\nu w \equiv 0$ and $[w,\nu] = \nabla_w \nu - \nabla_\nu w = \nabla_w \nu$. Since $\nabla_w \nu \perp \nu$ we can write $[w,\nu] = g \cdot w$ for some $C^1$-smooth function $g$ but 
$$0 = g \mathcal{L}_w |\nabla f|= \mathcal{L}_{g \cdot w} |\nabla f|= \mathcal{L}_{[w,\nu]} |\nabla f| = \mathcal{L}_w \mathcal{L}_\nu |\nabla f|- \mathcal{L}_\nu \mathcal{L}_w |\nabla f| = \mathcal{L}_w \mathcal{L}_\nu |\nabla f|$$
so $\mathcal{L}_w \mathcal{L}_\nu |\nabla f| \equiv \mathcal{L}_w |\nabla f| \equiv 0$ 
and $|\nabla f| , \mathcal{L}_\nu |\nabla f|$ are constant along regular level sets. This proves (3) and (4).\\\\
Let $c \in \mathbb{R}$ be some regular value. Then because $w , \nu$ are an eigenbasis and $w$ corresponds to the largest eigenvalue (in absolute value)  $|\nabla f||\nabla_w \nu|=|<H_f w ,w>|$ is greater then $|\mathcal{L}_\nu| = |<H_f \nu ,\nu>|$. Integration on the level set $f^{-1}(c)$ gives
$$ \int_{f^{-1}(c)} |L_\nu |\nabla f|| \: ds \leq \int_{f^{-1}(c)} |\nabla f||\nabla_w \nu| \: ds$$
and since $|\nabla f| , \mathcal{L}_\nu |\nabla f|$ are constant on the level set
$$ |L_\nu |\nabla f|| \cdot Length(f^{-1}(c))= |L_\nu |\nabla f|| \int_{f^{-1}(c)} 1 \: ds \leq |\nabla f| \int_{f^{-1}(c)}|\nabla_w \nu| \: ds = 2 \pi |\nabla f|$$
this proves (5). \\\\
Property (6) would remain as a debt until after Theorem \hyperlink{theorem:3.7}{3.7}. So we will not use it during this section at all. Given that we completed the proof of one direction. \\\\
On the other direction, we assume that for every regular value $c \in \mathbb{R}$ properties (1) to (6) hold and show that $\underset{M}{\max}|f| = \frac{1}{2\pi}\int_{M}||H_f||_{op}d\sigma$.\\\\
Note that from property (6) we know that $f \not\equiv 0$. From properties (1) and (2):
$$ B(1,f) = \int^{\infty}_{-\infty}\beta(c,f)dc = \int^{\underset{M}{\sup} f}_{\underset{M}{\inf} f} 1 dc = \underset{\overline{M}}{\max} f - \underset{\overline{M}}{\min} f = \underset{M}{\max} |f|.$$
Define $\mathcal{R}$ as the set of regular values of $f|_M$
and $U=f^{-1}(\mathcal{R})$. By Sard lemma $\mathcal{R}$ is open so $U$ is open as well. Denote $w = \dot{\gamma}(t)$ the tangent vector and $\nu = \pm \frac{\nabla f}{|\nabla f|}$ the outward unit normal.
The vectors $w,\nu$ form an orthonormal frame on $U$ so $\nabla_w \nu \perp \nu$. Using properties (3) and (4): 
$$H_f \nu = \nabla_{\nu}\nabla f = \nabla_{\nu}(|\nabla f| \nu) = |\nabla f| \nabla_{\nu}\nu + \mathcal{L}_\nu |\nabla f| \nu = \mathcal{L}_\nu |\nabla f| \nu$$ 
$$H_f w = \nabla_w\nabla f = \nabla_w(|\nabla f| \nu) = |\nabla f| \nabla_w \nu + \mathcal{L}_w |\nabla f| \nu = |\nabla f| \nabla_w \nu $$
so $w,\nu$ are an eigenbasis to $H_f$ on $U$ with corresponding eigenvalues $\lambda_w , \lambda_\nu$. Now, observe that $|\lambda_w| = |\nabla f| |\nabla_w \nu|$ and $|\lambda_\nu| = |\mathcal{L}_\nu |\nabla f||$ and from property (5) $|\lambda_w| \geq |\lambda_\nu|$ hence
$||H_f||_{op} = |\lambda_w| = |\nabla f| |\nabla_w \nu|$.\\\\
Using property (1) and Fenchel-Borsuk theorem again
\begin{align*}
L(c) = \frac{1}{2\pi}\int_{f^{-1}(c)}|\nabla_w \nu(s)|ds = \frac{1}{2\pi}\int_{f^{-1}(c)}|\kappa(s)|ds = 1 = \beta(c,f)
\end{align*}
and finally by utilizing property (6) we can use the smooth co-area formula 
\begin{align*}
& \frac{1}{2\pi}\int_{M}||H_f||_{op} d\sigma = 
\frac{1}{2\pi}\overset{\infty}{\underset{-\infty}{\int}} \left(\int_{f^{-1}(c)}\frac{||H_f||_{op}}{|\nabla f|}ds \right) dc = \\ &
\frac{1}{2\pi}\overset{\infty}{\underset{-\infty}{\int}} \left(\int_{f^{-1}(c)}\frac{|\nabla f| |\nabla_w \nu|}{|\nabla f|}ds \right) dc  = \frac{1}{2\pi}\overset{\infty}{\underset{-\infty}{\int}} \left(\int_{f^{-1}(c)} |\nabla_w \nu| ds \right) dc = \\ & \frac{1}{2\pi}\overset{\infty}{\underset{-\infty}{\int}} \beta(c,f) dc = B(1,f) = \underset{M}{\max} |f|
\end{align*}
and this completes the proof altogether.
\end{proof}
\end{Theorem2}
We will need a slightly stronger version of property $(1)$ in the following sections. 

Let $\gamma$ be a simple closed curve. Jordan curve theorem states that $\gamma$ divides $\mathbb{R}^2$ into two disjoint domains. Denote the domain bounded by $\gamma$ as \textbf{\textsl{int}}$(\gamma)$.
We denote the other (unbounded) domain by \textbf{\textsl{ext}}$(\gamma)$. \hypertarget{lemma:3.2}{} \\
\begin{Theorem16} Let $c \in \mathbb{R}$ be a regular value of $f$. If $f$ attains equality in (\hyperlink{eq:ineq}{1}) then: \\
(1*) $f^{-1}(c)$ is a convex plane curve. Denote by $\Omega$ the interior of $f^{-1}(c)$. Then $\Omega$ is a convex domain. Moreover,
$$\Omega = \left\{p \in M| f(p) > c\right\}.$$
\begin{proof}
We have already seen that $f^{-1}(c)$ is a convex plane curve, hence $\Omega$ is a convex domain. So all we need to show is that $\Omega = \left\{f > c\right\}$. Note that if the curve $f^{-1}(c)$ circles a hole (as in the figure below) then $\Omega$ is not even a subset of $\overline{M}$. Therefore we divide the proof into three parts: \\
\textbf{a.)} $\Omega \cap \overline{M} \subseteq  \left\{f > c\right\}$. \\
\textbf{b.)} $\Omega \subseteq \overline{M}$. \\
\textbf{c.)} $\left\{f > c\right\} \subseteq \Omega$. \\\\
\textbf{a.)} We need to show that \textbf{\textsl{int}}$\left\{f=c\right\} \cap \overline{M} \subseteq  \left\{f>c\right\}$. Assume that is not the case. Then there exists $x \in$ \textbf{\textsl{int}}$\left\{f=c\right\} \cap \overline{M}$ such that $x \notin \left\{f > c\right\}$. That is, $b:=f(x) < c$. 

Without loss of generality $x$ is a regular point.\\
Explanation: note that $\overline{M}$ and \textbf{\textsl{int}}$\left\{f=c\right\} \cup f^{-1}(c)$ are compact connected sets. Therefore their intersection is also compact and connected. Since $\mathbb{R}^2$ is locally connected the intersection is pathwise connected as well, so we can pick a curve $\delta_1:[0,1] \rightarrow \overline{M}$ such that:
$$\delta_1(0)=x \: ,\: \: \delta_1([0,1)) \in \mathbf{int} \left\{f=c\right\} 
\; \; \mathrm{and} \; \; \delta_1(1) \in \left\{f=c\right\}.$$
By definition $b,c \in (f \circ \delta_1)([0,1])$ and $f$ continuous so the segment $[b,c]$ is also in the image. Sard lemma guarantees that we can find a regular value $b' \in [b,c]$, but $b' \in [b,c]\subseteq (f \circ \delta_1)([0,1])$ so there exists some $t \in [0,1)$ such that 
$f(\delta_1(t))=b'$. By construction $x':=\delta_1(t) \in$ \textbf{\textsl{int}}$\left\{f=c\right\}$
and $b'=f(x')$ is a regular value. 
\begin{figure}[H]
	\centering
	\includegraphics[width=0.6\textwidth]{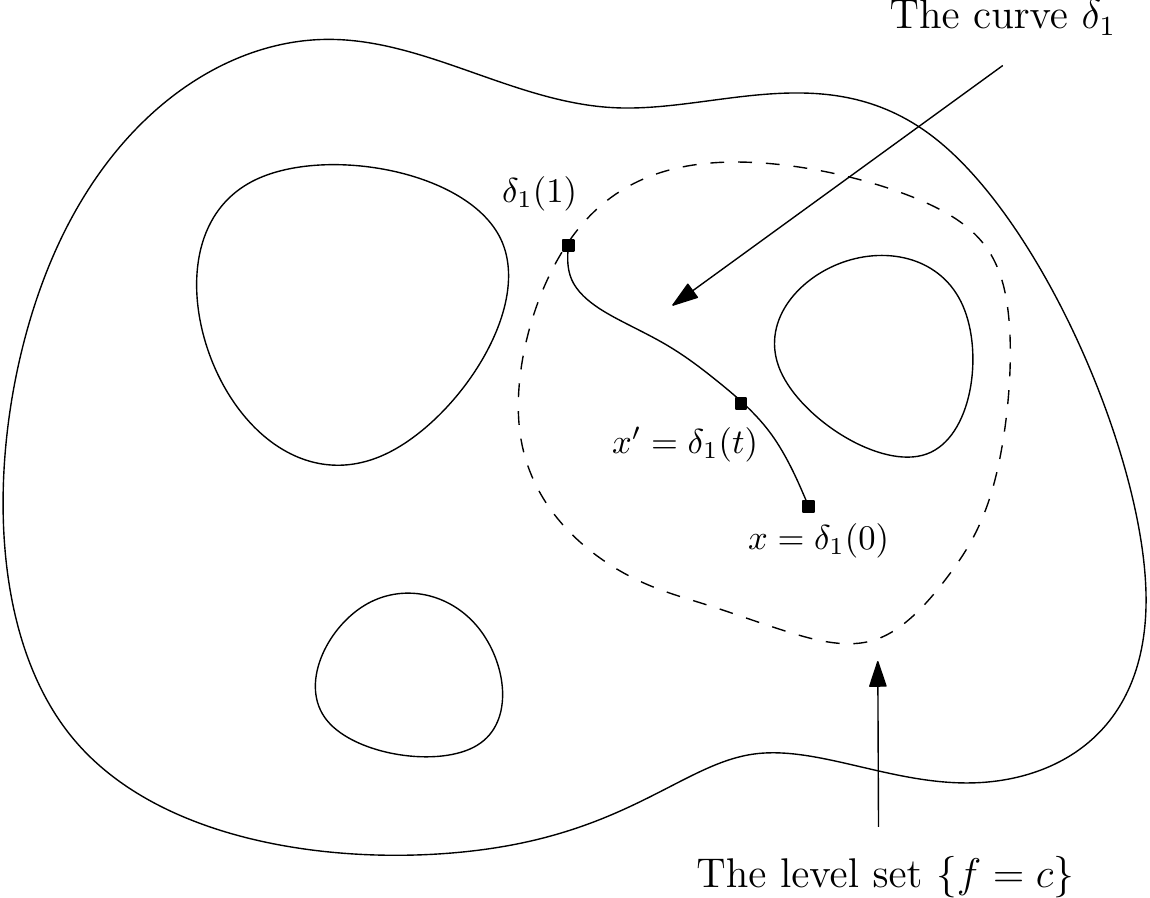}
	\caption{There exists $t \in [0,1]$ such that $\delta_1(t)$ is a regular point for $f$.}
\end{figure} 
From similar reasons, we can pick a curve $\delta_2:[0,1] \rightarrow \overline{M}$ such that:
$$\delta_2(0) \in \left\{f=c\right\} \: ,\: \: \delta_2((0,1]) \in \mathbf{ext} \left\{f=c\right\}
\; \; \mathrm{and} \; \; \delta_2(1) \in \partial M.$$
Note that because $c,0 \in (f \circ \delta_2)([0,1])$ the segment $[0,c] \subseteq (f \circ \delta_2)([0,1])$, so there exists $0<t<1$ such that $f(\delta_2(t))=b$. But from property (1) in proposition \hyperlink{theorem:3.1}{3.1} we know that the level set $f^{-1}(b)$ is connected.
Hence $f^{-1}(b) \subseteq \textbf{\textsl{ext}}\left\{f=c\right\}$. But $b$ was chosen in such a way that $b \in f(\textbf{\textsl{int}}\left\{f=c\right\})$ and contradiction follows. Therefore 
we have proven that \textbf{\textsl{int}}$\left\{f=c\right\} \cap \overline{M} \subseteq \left\{f>c\right\}$.\\\\
\textbf{b.)} Next, we show that $\Omega =$ \textbf{\textsl{int}}$\left\{f=c\right\} \subseteq \overline{M}$. Assume that is not the case. Then $\Omega \cap \overline{M}^c \neq \phi$. We observe that $$\Omega = (\Omega \cap M) \cup (\Omega \cap \overline{M}^c) \cup (\Omega \cap \partial M)$$
and the sets $\Omega \cap M$ and $\Omega \cap \overline{M}^c$ are open in $\Omega$. It is easy to see that $\Omega \cap M \neq \phi$. Since $\Omega$ is connected, it follows that $\Omega \cap \partial M$ is not empty. 

Let $y \in \Omega \cap \partial M$. Then $y \in \partial M$ and $f$ vanishes on the boundary so $f(y)=0<c$ and $y \in \Omega$. But we have already proved in case \textbf{a.)} that is not possible. So $\Omega \subseteq \overline{M}$. Combined with the previous case, this shows that $\Omega \subseteq \left\{f>c\right\}$.\\\\
\textbf{c.)} In order to complete the proof and show equality, we only need to prove that $\left\{f>c\right\} \subseteq$ \textbf{\textsl{int}}$\left\{f=c\right\}$. Assume that is not case, then there exists 
some $x \notin$\textbf{\textsl{int}}$\left\{f=c\right\}$ such that $f(x)>c$. Take any 
$y \in$\textbf{\textsl{int}}$\left\{f=c\right\} \subseteq \left\{f>c\right\}$. $\overline{M}$ is pathwise connected, so there exists a curve $\delta_3:[0,1] \rightarrow \overline{M}$ joining $x$ and $y$. Because $x \in$\textbf{\textsl{ext}}$\left\{f=c\right\}$ and $y \in$\textbf{\textsl{int}}$\left\{f=c\right\}$
there exists some $0<t<1$ such that $f(\delta_3(t))=c$. Take any regular value $c'\in \mathbb
{R}$ such that $c<c'<\min(f(x),f(y))$. It is obvious the value $c'$ is taken at least twice by $f \circ \delta_3$ so $\beta(c',f) \geq 2$. This is a contradiction to property (1) in proposition \hyperlink{theorem:3.1}{3.1}. So in conclusion $\left\{f>c\right\}=\Omega$ as requested.
\end{proof}
\end{Theorem16}
\subsection{Uniform change of $f$ and $|\nabla f|$ along the normal flow}
Let $M$ be an open and bounded domain in $\mathbb{R}^2$. Let $f \in C^3_0(M)$ be a non constant solution to the equation
$$\underset{M}{\max}|f| = \frac{1}{2\pi}\int_{M}||H_f||_{op}d\sigma.$$
Using property (2) in proposition \hyperlink{theorem:3.1}{3.1} we can assume, without loss of generality, that $f$ is non-negative. So from here on $f \geq 0$. \\\\
Denote $U = \left\{ x \in M \left|\right. |\nabla f|_x \neq 0 \right\}$. Note that the complement $\overline{M} \backslash U$ can be written as the union
$$\overline{M} \backslash U = \left\{ x \in \partial M \right\} \cup \left\{ x \in M \left|\right. |\nabla f|_x = 0 \right\}.$$ 
Since both sets are obviously closed,  $\overline{M} \backslash U$ is closed as well. Moreover, $\overline{M} \backslash U \subseteq \overline{M}$ and $\overline{M}$ is bounded so $\overline{M} \backslash U$ is compact and $U$ is open. Define the field of outward unit normals on $U$ by the rule $$\nu(x):=-\frac{\nabla f(x)}{|\nabla f|(x)}$$ for every $x \in U$. Note that $\nu$ is well defined and $C^2$-smooth because $f$ is $C^3$-smooth and $|\nabla f| \neq 0$ on $U$. Let $c \in \mathbb{R}$ be a regular value for $f$. Let $x_0 \in f^{-1}(c)$. Denote by $\gamma_{x_0}$ the flow associated with the vector field $\nu$ and emanating from $x_0$. That is,  $\gamma_{x_0}:\left(-a_{x_0},b_{x_0}\right) \rightarrow U$ is the solution to the O.D.E.
$$\dot{\gamma}_{x_0}(t) = \nu(\gamma_{x_0}(t))$$
with initial condition: $\gamma_{x_0}(0)=x_0$ (where $\mathbb{I}_{x_0}=\left(-a_{x_0},b_{x_0}\right)$ denotes the maximal interval of definition for the solution $\gamma_{x_0}$). Define the functions $g_{x_0}(t)=f(\gamma_{x_0}(t))$ and $h_{x_0}(t)=|\nabla f|(\gamma_{x_0}(t))$. We note that since $\nu$ is a $C^2$-smooth vector field, the integral curve $\gamma_{x_0}(t)$ is also in $C^2$ (see \cite[Theorem 17.19]{[2]}). By the composition law the functions $g_{x_0}(t)$ and $h_{x_0}(t)$ are $C^2$ as well.

The purpose of this section is prove proposition \hyperlink{theorem:3.7}{3.7}: the segment $\mathbb{I}_{x_0}$ and the functions $g_{x_0}(t)$ and $h_{x_0}(t)$ depend only on the regular value $c=f(x_0)$ and not on the choice of initial point $x_0$. We will use this proposition in the next section to show that $M$ is an open disc.

We begin the proof of proposition \hyperlink{theorem:3.7}{3.7} by making the following observation:
Define $\mathcal{R}_{x_0}:=\left\{ t \in \left(-a_{x_0},b_{x_0}\right) \left|\right. \gamma_{x_0}(t) \right.$ is a regular point for $f$ $\left.\right\}$. Let $\mu$ denote the standard Lebesgue measure on $\mathbb{R}$. Then
\begin{Theorem15} The set $\mathcal{R}_{x_0}$ is open and dense in $\mathbb{I}_{x_0}$ and $\mu(\mathbb{I}_{x_0} \backslash \mathcal{R}_{x_0})=0$.
\begin{proof}
Let $\mathcal{R}$ be the set of regular values of $f$. Denote by $\mathcal{S}$ the set of singular values. 
We make the following observation
$$\frac{d}{dt}f(\gamma_{x_0})(t) = <\nabla f|_{\gamma_{x_0}(t)} , \nu|_{\gamma_{x_0}(t)} > = - |\nabla f|(\gamma_{x_0}(t)) < 0$$
so $g_{x_0}:(-a_{x_0},b_{x_0}) \rightarrow \mathbb{R}$ is a strictly decreasing. Also, the set $Im(g_{x_0}) \subseteq \mathbb{R}$ is open and connected. So $Im(g_{x_0})$ is an open interval and we can denote $Im(g_{x_0}) = (-c_{x_0},d_{x_0})$. From the above, it follows that there is a well defined $C^2$-smooth inverse function $g_{x_0}^{-1}:(-c_{x_0},d_{x_0}) \rightarrow (-a_{x_0},b_{x_0})$. \\\\
According to the definition:
\begin{center}
$\mathcal{R}_{x_0} = \left\{ t \in \left(-a_{x_0},b_{x_0}\right) \left|\right. g_{x_0}(t) \right.$ is a regular value $\left.\right\}=g_{x_0}^{-1}(\mathcal{R} \cap (-c_{x_0},d_{x_0}))$.
\end{center}
By Sard lemma $\mathcal{R}$ is open and dense, so $g_{x_0}^{-1}(\mathcal{R} \cap (-c_{x_0},d_{x_0}))$ is open and dense as well. So what is left is to show that
$\mu(g_{x_0}^{-1}(\mathcal{S} \cap (-c_{x_0},d_{x_0})))=0$:
 
Let $[\alpha,\beta] \subset (-c_{x_0},d_{x_0})$. Note that because $g_{x_0}^{-1}$ is continuously differentiable,
$g_{x_0}^{-1}|_{[\alpha,\beta]}$ is Lipschitz and thus absolutly continuous. By Sard lemma $\mu(\mathcal{S})=0$ so $\mu(\mathcal{S} \cap [\alpha,\beta])=0$ as well. $g_{x_0}^{-1}$ is absolutly continuous so $\mu(g_{x_0}^{-1}(\mathcal{S} \cap [\alpha,\beta]))=0$. Since that is true for every $[\alpha,\beta]$ in $(-c_{x_0},d_{x_0})$, $\mu(g_{x_0}^{-1}(\mathcal{S} \cap (-c_{x_0},d_{x_0})))=0$.\\\\
\end{proof}
\end{Theorem15}  
We need the following two lemmata on $\mathbb{I}_{x_0}$: \hypertarget{lemma:3.4}{}
\begin{Theorem12} The segment $\left(-a_{x_0},b_{x_0}\right)$ is finite and $\gamma_{x_0}:\left(-a_{x_0},b_{x_0}\right) \rightarrow \mathbb{R}$ is a straight line. 
\begin{proof} The set $\mathcal{R}_{x_0}$ is open, so we can write it as the disjoint union of a countable number of open intervals:
$$\mathcal{R}_{x_0} = \bigcup_{i=0}^{\infty} (\alpha_i,\beta_i).$$
Note that from property (3) in \hyperlink{theorem:3.1}{3.1}, in every such interval $(\alpha_i,\beta_i)$: 
$$\nabla_{\dot{\gamma_{x_0}}} \dot{\gamma_{x_0}} = \nabla_\nu \nu = 0.$$
So $\gamma_{x_0}|_{(\alpha_i,\beta_i)}$ is a geodesic line for the euclidean metric i.e. a straight line.
Note that this means that the function $\dot{\gamma_{x_0}}(t)$ is constant on every such interval $(\alpha_i,\beta_i)$. Thus $\mu(\dot{\gamma_{x_0}}(\alpha_i,\beta_i))=0$. 

We want to show that $\dot{\gamma_{x_0}}:\mathbb{I}_{x_0} \rightarrow \mathbb{R}$ is constant. Assume that is not the case. Then there exists $\alpha,\beta \in \mathbb{I}_{x_0}$ such that $\dot{\gamma_{x_0}}(\alpha) < \dot{\gamma_{x_0}}(\beta)$. The curve $\gamma_{x_0}(t)$ is $C^2$ so $\dot{\gamma_{x_0}}|_{[\alpha,\beta]}$ is $C^1$-smooth and thus absolutly continuous. From Luzin property:
\begin{align*}
\mu(\dot{\gamma_{x_0}}([\alpha,\beta])) &= \mu(\dot{\gamma_{x_0}}([\alpha,\beta] \cap \mathcal{R}_{x_0}))
\cup \mu(\dot{\gamma_{x_0}}([\alpha,\beta] \cap \mathcal{R}_{x_0}^c)) \\ &=  
\sum_{i=0}^{\infty} \mu(\dot{\gamma_{x_0}}([\alpha,\beta] \cap (\alpha_i,\beta_i))) + \mu(\dot{\gamma_{x_0}}([\alpha,\beta] \cap \mathcal{R}_{x_0}^c)) \\ &= \sum_{i=0}^{\infty} 0 + 0 = 0.
\end{align*}
But $\dot{\gamma_{x_0}}$ continuous, so the image $\dot{\gamma_{x_0}}([\alpha,\beta])$ contains the inteval 
$[\dot{\gamma_{x_0}}(\alpha),\dot{\gamma_{x_0}}(\beta)]$. Obviously 
$$\mu(\dot{\gamma_{x_0}}([\alpha,\beta])) \geq \dot{\gamma_{x_0}}(\beta) - \dot{\gamma_{x_0}}(\alpha) >0 $$
which is a contradiction. So $\dot{\gamma_{x_0}}$ is constant and $\gamma_{x_0}:\mathbb{I}_{x_0} \rightarrow \mathbb{R}$ is a straight line. 

Because $M$ is a bounded domain every straight line contained in $M$ has finite length. In particuler
$$\int^{b_{x_0}}_{-a_{x_0}} 1 dt = \int^{b_{x_0}}_{-a_{x_0}}|\nu|dt = Length(\gamma_{x_0}) < \infty$$
so $a_{x_0},b_{x_0}$ are both finite.
\end{proof}
\end{Theorem12}
\hypertarget{lemma:3.5}{}
\begin{Theorem13} For $a_{x_0}$, the limit $\underset{t \rightarrow -a_{x_0}^+}{\lim} \gamma_{x_0}(t)$ exists and either 
$$f(\underset{t \rightarrow -a_{x_0}^+}{\lim} \gamma_{x_0}(t))=0 \: \mathrm{or} \: |\nabla f|(\underset{t \rightarrow -a_{x_0}^+}{\lim} \gamma_{x_0}(t))=0.$$
Moreover, the same is true for $b_{x_0}$ as well.
\begin{proof}
We have shown that $\gamma_{x_0}$ is a finite straight line so obviously $\underset{t \rightarrow -a_{x_0}^+}{\lim} \gamma_{x_0}(t)$ exists. If $\underset{t \rightarrow -a_{x_0}^+}{\lim} \gamma_{x_0}(t) \in U$ then we could easily extend the solution $\gamma_{x_0}$ to some larger interval $\mathbb{I}_{x_0} \subset \mathbb{I}'_{x_0}$ (see for example \cite[p. 102]{[3]}). These is a contradiction to the definition of $\mathbb{I}_{x_0}$ as a maximal interval of definition. So $\underset{t \rightarrow -a_{x_0}^+}{\lim} \gamma_{x_0}(t) \in \partial U$. The same argument shows that $\underset{t \rightarrow b_{x_0}^-}{\lim} \gamma_{x_0}(t) \in \partial U$. But $\partial U \subseteq \overline{M} \backslash U$ because $\overline{M} \backslash U$ is closed. The complement $\overline{M} \backslash U$ is the union of two sets 
$$\overline{M} \backslash U = \left\{ x \in \partial M \right\} \cup \left\{ x \in M \left|\right. |\nabla f|_x = 0 \right\}.$$ 
Now the result follows from the fact that $f$ vanishes on the boundary.
\end{proof}
\end{Theorem13}
Also, we would make use of the following auxiliary definitions and results: 
Let $\Psi$ be the vector field $\frac{\nu}{|\nabla f|} = - \frac{\nabla f}{|\nabla f|^2}$ on $U$.
Denote by $\delta_{x_0}:\left(-a_{x_0}',b_{x_0}' \right) \rightarrow U$ the integral curve of $\Psi$ with respect to the initial condition $\delta_{x_0}(0)=x_0$ (where $\left(-a_{x_0}',b_{x_0}' \right)$ is the maximal interval of definition). Note that $\delta_{x_0}$ is $C^2$-smooth, because the vector field $\Psi$ is $C^2$-smooth.

We have the following simple connection between $\gamma_{x_0}$ and $\delta_{x_0}$ \hypertarget{lemma:3.6}{}
\begin{Theorem14} Let $t \in \mathbb{I}_{x_0}$. Denote $\alpha = c - f(\gamma_{x_0}(t))$. Then $\alpha \in \left(-a_{x_0}',b_{x_0}' \right)$ and $\delta_{x_0}(\alpha)=\gamma_{x_0}(t)$.
\begin{proof} First, note that $\delta_{x_0}$ has the same image as $\gamma_{x_0}$ so necessarly there exists a solution to the equation $\delta_{x_0}(\alpha)=\gamma_{x_0}(t)$. Denote such a solution by $\alpha_0$.
We observe that $(f \circ \delta_{x_0})' = <\nabla f , -\frac{\nabla f}{|\nabla f|^2}> = -1$. By the definition of $\alpha$ and $\alpha_0$:
\begin{align*}
\alpha = c - f(\gamma_{x_0}(t)) &=(f \circ \delta_{x_0})(0) - (f \circ \delta_{x_0})(\alpha_0) = -\int_0^{\alpha_0} (f \circ \delta_{x_0})'(\tau) d\tau \\ &= \int_0^{\alpha_0} 1 d\tau  = \alpha_0 
\end{align*}
as required.
\end{proof}
\end{Theorem14}
We can now prove the proposition. \hypertarget{theorem:3.7}{}
\begin{Theorem4} Let $x_0,y_0 \in f^{-1}(c)$. Then:\\
(1) For every $t \in \mathbb{I}_{x_0} \cap \mathbb{I}_{y_0}$:
$$ f(\gamma_{x_0}(t))=f(\gamma_{y_0}(t)) \: \: \mathrm{and} \: \: |\nabla f|(\gamma_{x_0}(t)) = |\nabla f|(\gamma_{y_0}(t)).$$
(2) $\mathbb{I}_{x_0}=\mathbb{I}_{y_0}$ and $\mathcal{R}_{x_0}=\mathcal{R}_{y_0}$.
\begin{proof} (1) Since $0 \in \mathbb{I}_{x_0} \cap \mathbb{I}_{y_0}$, it is obvious that 
$\mathbb{I}_{x_0} \cap \mathbb{I}_{y_0} \neq \phi$. Let $[\alpha,\beta] \subseteq \mathbb{I}_{x_0} \cap \mathbb{I}_{y_0}$ be a closed subinterval. Define a smooth function $p : \left(f(\gamma_{x_0}(\alpha)),f(\gamma_{x_0}(\beta)) \right) \rightarrow \mathbb{R}$ as the composition $p(t) = |\nabla f|(\delta_{x_0}(c - t))$. From lemma \hyperlink{lemma:3.6}{3.6} the function $p$ is well defined. As a composition of $C^2$-smooth functions $p$ itself is $C^2$-smooth. Note that if $t$ is a regular value of $f$ it follows from property (4) in Proposition \hyperlink{theorem:3.1}{3.1} that $p(t)$ is the common value of $|\nabla f|$ on $f^{-1}(t)$.

Define the set $\mathcal{R}:=\mathcal{R}_{x_0} \cap \mathcal{R}_{y_0}$. As the intersection of two dense sets, the set $\mathcal{R}$ is dense in $[\alpha,\beta]$ as well. For every $t \in \mathcal{R}$ (and thus from continuity for every $t \in [\alpha,\beta]$):
$$\dot{g}_{x_0}(t) = |\nabla f|(\gamma_{x_0}(t)) = p(f(\gamma_{x_0})(t)) = p(g_{x_0}(t))$$
$$\dot{g}_{y_0}(t) = |\nabla f|(\gamma_{y_0}(t)) = p(f(\gamma_{y_0})(t)) = p(g_{y_0}(t))$$
also $g_{x_0}(0)=g_{y_0}(0)=c$ implies $p(g_{x_0}(0))=p(g_{y_0}(0))$.  
Both $\dot{g}_{x_0},\dot{g}_{y_0}$ are solutions to the same O.D.E. with the same initial condition, therefore $f(\gamma_{x_0}(t))=g_{x_0}(t) = g_{y_0}(t)=f(\gamma_{y_0}(t))$ for all $t \in [\alpha,\beta]$. Thus (since the gradient is constant on regular level sets) for every $t \in \mathcal{R}$: $|\nabla f|(\gamma_{x_0}(t))=h_{x_0}(t) = h_{y_0}(t)=|\nabla f|(\gamma_{y_0}(t))$. Using the density of $\mathcal{R}$ in $[\alpha,\beta]$ once more, we can extend this result to every $t \in [\alpha,\beta]$ as well. But $[\alpha,\beta]$ was an arbitrary segment in $\mathbb{I}_{x_0} \cap \mathbb{I}_{y_0}$, so for every $t \in 
\mathbb{I}_{x_0} \cap \mathbb{I}_{y_0}$ the functions $g_{x_0}(t)=g_{y_0}(t)$ and $h_{x_0}(t)=g_{y_0}(t)$. 
Using the characterization above, it is obvious that $b_{x_0} = b_{y_0},a_{x_0} = a_{y_0}$. Therefore $g_{x_0}(t),h_{x_0}(t)$ are independent of the choice of $x_0$. \\\\
(2) Follows immediatly from (1) in \hyperlink{theorem:3.7}{3.7} and lemma \hyperlink{lemma:3.5}{3.5}.
\end{proof}
\end{Theorem4}
To conclude this subsection, we use the proposition above to make the following definition, \\\\
\textbf{Definition:} Let $c \in \mathbb{R}$ be a regular value. Define a segment $\mathbb{I}_c$, a set $\mathcal{R}_c$ and functions $g_c,h_c : \mathbb{I}_c \rightarrow \mathbb{R}$ in the following manner: take $x_0 \in f^{-1}(c)$. Then $\mathbb{I}_c = \mathbb{I}_{x_0},\mathcal{R}_{c} = \mathcal{R}_{x_0}$ and $g_c = g_{x_0} \: , \: h_c = h_{x_0}$. Proposition \hyperlink{theorem:3.7}{3.7} shows that $\mathbb{I}_c,\mathcal{R}_c$ and $g_c,h_c$ are well defined.
\subsection{$M$ is a open disc}
Let $M$ be an open and bounded domain. Let $f \in C^3_0(M)$ be a non constant solution to the equation
$$\underset{M}{\max}|f| = \frac{1}{2\pi}\int_{M}||H_f||_{op}d\sigma.$$
We assume that $f(x) \geq 0$ for all $x \in \overline{M}$. Let $c \in \mathbb{R}$ be a regular value of $f$. In the previous section we have defined a function 
$h_c : \mathbb{I}_c \rightarrow \mathbb{R}$. Using \hyperlink{theorem:3.1}{3.1} and \hyperlink{lemma:3.2}{3.2} we deduce the following estimates: \hypertarget{lemma:3.8}{}
\begin{Theorem17} (1) For every $t \in \mathcal{R}_c$ 
$$|\dot{h}_c(t)| \leq \frac{2\pi}{Length(f^{-1}(g_c(t)))} |h_c(t)|.$$ \\
(2) Let $t_1,t_2 \in \mathcal{R}_c$. If $t_1<t_2$ then 
$$\frac{2\pi}{Length(f^{-1}(g_c(t_2)))}<\frac{2\pi}{Length(f^{-1}(g_c(t_1)))}. $$
\begin{proof} (1) Take $x_0 \in f^{-1}(c)$. The curve $\gamma_{x_0}$ is an integral curve of the vector field $\nu$ so $$\mathcal{L}_\nu |\nabla f|(\gamma_{x_0}(t)) = \frac{d}{dt}|\nabla f|(\gamma_{x_0}(t)) = \dot{h}_{x_0}(t).$$
Let $t \in \mathcal{R}_c=\mathcal{R}_{x_0}$. Note that by the definition of $\mathcal{R}_{x_0}$ the point $\gamma_{x_0}(t)$ belongs to a regular level set of $f$. Using property (5) from \hyperlink{theorem:3.1}{3.1} we see that:
\begin{align*}
|\dot{h}_{x_0}(t)| &= |\mathcal{L}_\nu |\nabla f|(\gamma_{x_0}(t))| \leq
\frac{2\pi}{Length(f^{-1}(f(\gamma_{x_0}(t))))}|\nabla f|(\gamma_{x_0}(t)) \\ &=  \frac{2\pi}{Length(f^{-1}(g_{x_0}(t)))}|\nabla f|(\gamma_{x_0}(t))
= \frac{2\pi}{Length(f^{-1}(g_{x_0}(t)))}|h_{x_0}(t)| 
\end{align*}
Since $g_{c} = g_{x_0}$ and $h_{c} = h_{x_0}$ this is the required result. \\\\
(2) Take $x_0 \in f^{-1}(c)$. Let $t_1,t_2 \in \mathcal{R}$. Assume that $t_1<t_2$. Then $g(t_1),g(t_2)$ are both regular values. Note that 
$$\frac{d}{dt}f(\gamma_{x_0})(t) = <\nabla f|_{\gamma_{x_0}(t)} , \nu|_{\gamma_{x_0}(t)} > = - |\nabla f|(\gamma_{x_0}(t)) < 0$$
so $g_c = g_{x_0}$ is strictly decreasing and $g_c(t_1) < g_c(t_2)$. From property (1*) in \hyperlink{lemma:3.2}{3.2} the level sets $f^{-1}(g(t_1))$ and $f^{-1}(g(t_2))$ are both convex and 
$f^{-1}(g(t_1))$ encircles $f^{-1}(g(t_2))$. So obviously 
$$Length(f^{-1}(g(t_1)) < Length(f^{-1}(g(t_2)).$$
But this means that
$$ \frac{2\pi}{Length(f^{-1}(g(t_1)))} < \frac{2\pi}{Length(f^{-1}(g(t_2)))}$$
as required.
\end{proof}
\end{Theorem17} 
We also need the following elementary lemma from calculus: \hypertarget{lemma:3.9}{}
\begin{Theorem3} Let $\mathbb{I}$ be a half open segment such that $0 \in \mathbb{I}$ (i.e. $\mathbb{I} =(-a,0],-\infty \leq a <0$ or $\mathbb{I} = [0,b),0< b <\infty $). Let
$h : \mathbb{I} \rightarrow \mathbb{R}$ be a smooth never zero function and let $K>0$ be such that for $t \in \mathbb{I}$ the differential inequality $|\dot{h}(t)| \leq K  |h(t)|$ holds. Then $e^{-K |t|} \leq \frac{h(t)}{h(0)} \leq e^{K |t|}$ for every $t \in \mathbb{I}$.
\begin{proof}
Since $h(t) \neq 0$ for every $t \in \mathbb{I}$ we have the following inequality
$$|\dot{h}(t)| \leq K |h(t)|  \Rightarrow |\ln(|h|)'|=|\frac{\dot{h}}{h}|(t) \leq K $$
and 
$$|\ln(|h|)(t) - \ln(|h|)(0)| = |\int^{t}_{0} \ln(|h|)' ds| \leq |\int^{t}_{0} |\ln(|h|)'| ds| \leq
K |t|$$ 
so
$$ -K  |t| \leq \ln(|h|)(t) - \ln(|h|)(0) \leq K |t| $$
taking $e$ to the power of both sides
$$e^{-K |t|} \leq |\frac{h(t)}{h(0)}| \leq e^{K |t|}$$
$h(t)$ is never zero therefore it preserves sign and $|\frac{h(t)}{h(0)}| = \frac{h(t)}{h(0)}$. This yields the required result. 
\end{proof}
\end{Theorem3} 
Now we are ready to prove the theorem.\\\\
\textbf{Proof of theorem \hyperlink{theorem:1.1}{1.1}:} As discussed before, without loss of generality, we can assume that $f \geq 0$. We will divide the proof into four parts.\\
\textbf{a.)} There are no singular points in $\left\{ x \in M \left|\right. 0 < f(x) < \underset{M}{\max}f \right\}$.\\
\textbf{b.)} $f^{-1}(\underset{M}{\max} f)$ contains exactly one point.\\
\textbf{c.)} $f^{-1}(0) = \partial M$ and $f^{-1}\left((0,\underset{M}{\max} f]\right)=M$. \\
\textbf{d.)} Denote $f^{-1}(\underset{M}{\max} f)=\left\{ p \right\}$ then $M$ is an open disc around $p$.\\\\
\textbf{a.)} Let $c \in \mathbb{R}$ be a regular value of $f$. Denote $(-a_{c},b_{c}) = \mathbb{I}_{c}$. 
Let $x_0 \in f^{-1}(c)$ and denote $q_{x_0}=\underset{t \rightarrow b_{c}^-}{\lim} \gamma_{x_0}(t)$. 
We want to show that $q_{x_0} \in \partial M$. Assume that is not the case, that is $q_{x_0} \notin \partial M$. Since $\overline{M}$ is closed the limit point $q_{x_0}$ is in $\overline{M}$. So $q_{x_0} \in \overline{M} \backslash \partial M = M$. Using lemma \hyperlink{lemma:3.5}{3.5} we deduce that $f(q_{x_0})=0$ or $|\nabla f|(q_{x_0})=0$. Note that $f \geq 0$ and $q_{x_0} \in M$, so if $f(q_{x_0})=0$ then $q_{x_0}$ is a global minima. So anyway $|\nabla f|(q_{x_0})=0$. Denote $K_{c} = \frac{2\pi}{Length(f^{-1}(c))}$. From lemma \hyperlink{lemma:3.8}{3.8} for every $0< t \in \mathcal{R}_{c}$ we have an estimate: 
$$|\dot{h}_c(t)| \leq \frac{2\pi}{Length(f^{-1}(g_c(t)))}|h_c(t)| \leq \frac{2\pi}{Length(f^{-1}(g_c(0)))}|h_c(t)| = K_c|h_c(t)|.$$
Note that the set $\mathcal{R}_{c}$ is dense in $\mathbb{I}_{c}$, so the estimate above actually holds for every $0< t < b_{c}$. Now we can use lemma \hyperlink{lemma:3.9}{3.9} to bound $h_{c}(t)$ from below:
$$\forall t \: \: \mathrm{s.t.} \: \: 0< t < b_{c} \: \: : \: \: \: 0< e^{-K_{c} b}h_{c}(0) \leq e^{-K_{c} t}h_{c}(0) \leq h_{c}(t).$$ 
But $$|\nabla f|(q_{x_0}) = |\nabla f|(\underset{t \rightarrow b_{c}^-}{\lim} \gamma_{x_0}(t)) = 
\underset{t \rightarrow b_{c}^-}{\lim} |\nabla f|(\gamma_{x_0}(t)) \stackrel{by \: def.}{=}
\underset{t \rightarrow b_{c}^-}{\lim} h_{x_0}(t) = \underset{t \rightarrow b_{c}^-}{\lim} h_{c}(t).$$
This is a contradiction, because we have shown 
$$\underset{t \rightarrow b_{c}^-}{\lim} h_{c}(t)>e^{-K_{c} b}h_{c}(0) >0$$
while necessarily $|\nabla f|(q_{x_0})=0$. So $\underset{t \rightarrow b_{c}^-}{\lim} \gamma_{x_0}(t) \in \partial M$ for every $x_0$ in a regular level set $f^{-1}(c)$. 

Let $c_1<c_2$ be two regular values of $f$. Then $\underset{t \rightarrow b_{c_2}^-}{\lim} g_{c_2}(t) = 0$ and $g_{c_2}(0)=c_2$ so there exists an intermidate value $0<t_1<b_{c_2}$ such that $g_{c_2}(t_1)=c_1$. We define the following homotopy between the closed curves $f^{-1}(c_1)$ and $f^{-1}(c_2)$:
$$H(s,t) : f^{-1}(c_2) \times [0,t_1] \rightarrow M$$
is defined by the rule $H(s,t) = \gamma_{s}(t)$. Let $0<t<t_1$. The common value of the gradient $|\nabla f|$ on the level set $f^{-1}(t)$ is $h_{c_2}(t)$. Observe that since $t_1<b_{c_2}$ the gradient $h_{c_2}(t) \neq 0$. Therefore $f^{-1}(t)$ is a regular level set for every $0<t<t_1$ and $H(s,t)$ is a regular homotopy. We want to show that $$Im H := \left\{H(s,t)\left|\right. s \in f^{-1}(c_2) \: \: \mathrm{and} \: \: t \in [0,t_1] \right\} = f^{-1}[c_1,c_2].$$
It is obvious that $Im H\subseteq f^{-1}[c_1,c_2]$. Assume that equality does not hold, then there exists $p \in M$ such that $c_1 < f(p) < c_2$ and $p \notin Im H$. Note that from property (1*) in \hyperlink{lemma:3.2}{3.2} the winding number of $f^{-1}(c_1)$ and $f^{-1}(c_2)$ around $p$ are different: $W(f^{-1}(c_1),p)=0$ and $W(f^{-1}(c_2),p)= \pm 1$. But $H(s,t)$ is a regular homotopy between $f^{-1}(c_1)$ and $f^{-1}(c_2)$ that does not pass through $p$ and as such preserves the winding number! we arrived at a contradiction, so $Im H = f^{-1}[c_1,c_2]$.

Note that $H$ is a regular homotopy so $Im H$ does not contain singular points. So for any two regular values $c_1 < c_2$ the set $f^{-1}[c_1,c_2]$ does not contain singular points. From Sard lemma the regular values are dense in $[0,\underset{M}{\max} f]$ and the result follows. \\\\
\textbf{b.)} We want to show that $|f^{-1}(\underset{M}{\max} f)|=1$. Assume that is not the case. Since $f$ is continuous, non-negative and $\overline{M}$ is closed the set $f^{-1}(\underset{M}{\max} f)$ is obviously not empty. So according to the assumption there exists $p \neq q \in M$ such that $f(p)=f(q)=\underset{M}{\max}f$. 

Let $c \in \mathbb{R}$ be a regular value of $f$. Because $\underset{M}{\max}f > c$ the points $p,q \in \textbf{\textsl{int}}(f^{-1}(c))$. Denote that straight line that connects $p$ and $q$ by $L$; also denote the length of $L$ by $\ell$. Note that since $f^{-1}(c)$ is a convex set the line $L$ is contained in $\textbf{\textsl{int}}(f^{-1}(c))$ and therefore $Length(f^{-1}(c)) \geq \ell$. This is a crucial part of the proof. The lower bounds on lengths of regular level sets allows us to use \hyperref[lm:3.8]{3.8} and \hyperref[lm:3.9]{3.9} to bound $h_c(t) , -a_c<t  \leq 0$ from below:
$$0< e^{(-\ell) \cdot (-a_c)} \leq e^{-\ell t} \leq \frac{h_c(t)}{h_c(0)} \Rightarrow 0 < e^{\ell \cdot a_c}h_c(0) \leq h_c(t).$$

Let $x_0 \in f^{-1}(c)$. Denote $p_{x_0} = \underset{t \rightarrow -a_c^+}{\lim} \gamma_{x_0}(t)$. Then
$|\nabla f|(p_{x_0}) \geq e^{\ell \cdot a_c}h_c(0) >0$ and therefore $\underset{t \rightarrow -a_c^+}{\lim} g_c(t) = f(p_{x_0})=0$. But $g_c$ is strictly decreasing so $0< c = g_{c}(0) < g_{c}(t)$ for every $-a_c<t \leq 0$ and contradiction follows. \\\\
\textbf{c.)} $f$ vanishes on the boundary, so obviously $\partial M \subset f^{-1}(0)$. So if we prove that $f^{-1}(0) \subset \partial M$ we are done. Assume that is not true. Then there exists 
$y_0 \in M = \overline{M} \backslash \partial M$ such that $f(y_0)=0$. Since $f$ is non-negative, the point $y_0$ is a local minima for $f$ and $|\nabla f|(y_0)=0$. But according to \textbf{a.)} and \textbf{b.)} this means that $0 = f(y_0) = \underset{M}{\max} f$. But $f$ is not constant so $f^{-1}(0)=\partial M$. \\\\
\textbf{d.)} Let $c \in \mathbb{R}$ be a regular value for $f$ and let $x_0 \in f^{-1}(c)$. Denote $p_{x_0} = \underset{t \rightarrow -a_c^+}{\lim}(\gamma_{x_0}(t))$. From lemma \hyperlink{lemma:3.5}{3.5} either $f(p_{x_0})=0$ or $|\nabla f|(p_{x_0})=0$. But $g_c(t)$ is strictly monotone decreasing so $\underset{t \rightarrow -a_c^+}{\lim} g_c(t) \geq g_c(0) = c >0$. Therefore $|\nabla f|(p_{x_0})=0$ and $\underset{t \rightarrow -a_c^+}{\lim}(\gamma_{x_0}(t)) \notin \partial M$. From \textbf{a.)} and \textbf{b.)} above it follows that $p_{x_0} \in f^{-1}(\underset{M}{\max}f)$ and $f^{-1}(\underset{M}{\max}f)=\left\{ p \right\}$, so $p = p_{x_0}$. Note that the distance from $x_0$ to $p$ depends only on $c = f(p)$ :
$$d(x_0,p) = \int_{-a_c}^{0} |\dot(\gamma_{x_0})(t)| dt = \int_{-a_c}^{0} 1 dt = a_c.$$ 
Denote by $\Omega_c$ the convex domain bounded by $f^{-1}(c)$. Every level set $f^{-1}(c)$ is a circle, so $\Omega_c$ is an open disc around $p$.

We note that the set $M = f^{-1} \left( (0,\underset{M}{\max} f] \right)$ can be written as the union
$$f^{-1} \left( (0,\underset{M}{\max} f]  \right) = \underset{regular \: c \in \mathbb{R}}{\bigcup} 
f^{-1}\left( (c,\underset{M}{\max} f] \right) = \underset{regular \: c \in \mathbb{R}}{\bigcup} \Omega_c.$$
By lemma \hyperref[lm:3.2]{3.2} this is the monotone union of open discs centered around $p$, so $f^{-1} \left( (0,\underset{M}{\max} f]  \right)$ is an open disc around $p$ as well. $\square$ \\\\ 
Also note that now we can repay a debt from proposition \hyperlink{theorem:3.1}{3.1}. In the last theorem we have shown that 
$$\left\{ x \in M \left|\right. {|\nabla f|}_x=0 \right\}=\left\{p\right\}$$ so $\mu(\left\{ x \in M \left|\right. {|\nabla f|}_x=0 \right\})=0$ which is exactly property (6).

\section{The equality case : description of the solutions}
\subsection{The sets $\mathcal{S}$ and $\mathcal{G}$}
As discussed in remark \hyperlink{theorem:1.2}{1.2} we can assume, without loss of generality, that $M=D(0,1)$.
Denote by $\mathcal{S} \subset C^3_0(D(0,1))$ the family of non constant solutions to the equality $$\underset{D(0,1)}{\max}|f| = \frac{1}{2\pi}\int_{D(0,1)}||H_f||d\sigma$$
with the normalization $f(0,0)=1$. Now define the following family of functions:
$$\mathcal{G}=\left\{g:[0,1] \rightarrow \mathbb{R}\left|\right. \: \: \mathrm{g \: \: is  \: \: C^2-smooth \: \: and}\: \: \forall t \in [0,1): -\frac{1}{t} \leq g'(t) \leq 0 \right\}$$
with the normalization condition $\int_0^{1} e^{g(t)}dt=1$.

In this section, we prove theorem \hyperlink{theorem:1.3}{1.3}: let $g \in \mathcal{G}$ and denote $h(t)=1 - \int^{t}_{0} e^{g(\tau)}d\tau$. Now take $f_g(x,y) = h(x^2+y^2)$. Then the correspondence 
$$\eta : \mathcal{G} \rightarrow \mathcal{S} \: \: , \: \: \eta(f) = f_g$$
is one to one from $\mathcal{G}$ onto $\mathcal{S}$. 
\hypertarget{lemma:4.1}{} 
\begin{Theorem5} Let $g \in \mathcal{G}$. Denote $f=f_g=\eta(g)$. Then:\\
(1) $|\nabla f| = -2\sqrt{x^2+y^2} \: h'(x^2+y^2)$.\\
(2) $|\mathcal{L}_\nu |\nabla f|| = 2|h'(x^2+y^2) + 2(x^2+y^2)h''(x^2+y^2)|$.\\
\begin{proof}
(1) $$f_x = \frac{\partial}{ \partial x} h(x^2+y^2) = 2x h'(x^2+y^2)$$
$$f_y = \frac{\partial}{ \partial y} h(x^2+y^2) = 2y h'(x^2+y^2)$$
$$|\nabla f| = \sqrt{f_x^2+f_y^2} = 2\sqrt{x^2+y^2} \: |h'(x^2+y^2)|= -2\sqrt{x^2+y^2} \: h'(x^2+y^2)$$\\
(2) $$\frac{\partial}{ \partial x} |\nabla f| = -2\left[\frac{x}{\sqrt{x^2+y^2}} h'(x^2+y^2) + 
\sqrt{x^2+y^2} h''(x^2+y^2)2x \right]$$
$$\frac{\partial}{ \partial y} |\nabla f| = -2\left[\frac{y}{\sqrt{x^2+y^2}} h'(x^2+y^2) + 
\sqrt{x^2+y^2} h''(x^2+y^2)2y \right]$$
The normal at the point $(x,y)$ is $\nu=-\frac{(x,y)}{\sqrt{x^2+y^2}}$ so the Lie derivative is given by the expression:
\begin{align*}
|\mathcal{L}_\nu |\nabla f|| &= |\frac{x}{\sqrt{x^2+y^2}} \cdot \frac{\partial}{\partial x}|\nabla f| + \frac{y}{\sqrt{x^2+y^2}} \cdot \frac{\partial}{\partial y}|\nabla f|| \\ &= 2|h'(x^2+y^2) + 2(x^2+y^2)h''(x^2+y^2)|.
\end{align*}\\
\end{proof}
\end{Theorem5}

\textbf{\\ Proof of theorem \hyperlink{theorem:1.3}{1.3}.} (1) Let $g \in \mathcal{G}$. Define the function $h:[0,1] \rightarrow \mathbb{R}$ as $h(t)=1 - \int^{t}_{0} e^{g(\tau)}d\tau$. $h$ is obviously continuous on $[0,1]$. Because $g \in C^2([0,1))$, $h$ is three times continuously differentiable in $[0,1)$.
So $$\eta(g) = f_g(x,y)=h(x^2+y^2)$$ is continuous in $\overline{D}(0,1)$, $C^3$ smooth in $D(0,1)$ and vanishes on the boundary. Moreover, 
$$f_g(0,0) = h(0) = 1 - \int^{0}_{0} e^{g(\tau)}d\tau = 1.$$
We want to show that $\eta(g) = f_g \in \mathcal{S}$. 

$\eta(g)$ clearly has properties (1)-(4) and (6) from theorem \hyperlink{theorem:3.1}{3.1}. So if we prove that $\eta(g)$ has property (5) we are done. 

We can express property (5) as a condition on $h$ (lemma \hyperlink{lemma:4.1}{4.1}):
$$2|h'(x^2+y^2) + 2(x^2+y^2)h''(x^2+y^2)| \leq -2h'(x^2+y^2).$$
After developing the left hand side and denoting $t=x^2+y^2$, 
$$h'(t) \leq h'(t) + 2t h''(t) \leq -h'(t) \: \:,\: \: 0<t<1.$$ 
Substituting $h'(t)=-e^{g(t)} \: , h''(t)=-g'(t) \cdot e^{g(t)}$ into the inequality yields an inequality in $g$:
$$|e^{g(t)} + 2tg'(t) \cdot e^{g(t)}| \leq e^{g(t)}$$
$e^{g(t)}$ is non vanishing, so we can divide both sides by $e^{g(t)}$ and get
$$-1 \leq 1 + 2tg'(t) \leq 1.$$
But $g \in \mathcal{G}$ so $-\frac{1}{t} \leq g'(t) \leq 0$ and the inequality holds. \\\\
(2) Let $f \in \mathcal{S}$. We want to show that there exists a unique $g \in \mathcal{G}$ such that $f_g=f$.\\\\
\textit{Existence:} The level sets of $f$ are circles around the origin, so $f$ is rotation invariant. 
It is easy to deduce (using the taylor expansion of $f$, for example) that there exists a continuous $h:[0,1] \rightarrow \mathbb{R}$, $C^3$-smooth on $[0,1)$, such that $f(x,y)=h(x^2+y^2)$.

Observe that from property (2) in theorem \hyperlink{theorem:3.1}{3.1} we know that $f$ is either non-positive or non-negative. Because $f(0,0)=1$ we conclude that $f$ is non-negative. So for every $t \in [0,1]$ the function $h$ is non-negative and for $t \in (0,1)$ the derivative $h'(t)$ is negative.

Using lemma \hyperlink{lemma:4.1}{4.1} again, we express property (5) in \hyperlink{theorem:3.1}{3.1} as a differential inequality in $h$:
$$h'(t) \leq h'(t) + 2t h''(t) \leq -h'(t) \: \: , \: \: 0<t<1.$$
Since $h'(t)<0$ for all $0 < t < 1$ the function $g(t):=\ln(-h'(t))$ is well defined. Note that $g:[0,1) \rightarrow \mathbb{R}$ is $C^2$-smooth, $\dot{h}(t)=-e^{g(t)}$ and $h''(t)=-g'(t)e^{g(t)}$. We can express the differential inequality above in terms of the function $g$:
$$-e^{g(t)} \leq -e^{g(t)} - 2tg'(t)e^{g(t)} \leq e^{g(t)} \: \: , \: \: 0<t<1.$$
Since $e^{g(t)}$ is never zero, we can divide both sides of the inequality by $e^{g(t)}$
$$-1 \leq -1 - 2tg'(t) \leq 1 \: \: , \: \: 0<t<1.$$
That is, $-\frac{1}{t} \leq g'(t) \leq 0$ for all $t \in (0,1)$. We also note that
$$\int^{1}_{0} e^{g(t)}dt = -\int^{1}_{0} \dot{h}(t) dt = h(0)-h(1) = f(0,0) - f(1,0)= 1.$$ 
So $g \in \mathcal{G}$ and
$$f_g(x,y) = 1 - \int^{x^2+y^2}_{0} e^{g(t)}dt = 1 + \int^{x^2+y^2}_{0} \dot{h}(t) dt = 1 + h(x^2+y^2) - h(0) = f(x,y)$$
so $g$ is the required solution and $\eta$ is onto. \\\\
\textit{Uniqueness:} Let $g_1,g_2$ be two functions in $\mathcal{G}$ such that $\eta(g_1)=\eta(g_2)=f$. Denote $h_1(t)=1 - \int^{t}_{0} e^{g_1(t)}dt$ and $h_2(t)=1 - \int^{t}_{0} e^{g_2(t)}dt $. 
By the assumption, $h_1(x^2+y^2) = \eta(g_1) \equiv \eta(g_2) = h_2(x^2+y^2)$ so $h_1 \equiv h_2$.
Taking the derivative of both sides of the equality we find that $e^{g_1(t)} = h_1'(t) = h_2'(t) = e^{g_2(t)}$ for all $t \in [0,1]$. So $g_1 \equiv g_2$ as requested. $\square$ 
\subsection{Examples}
Using the results of the previous section, we will now present some examples of functions from $\mathcal{G}$ and 
the corresponding solutions in $\mathcal{S}$:\\ \hypertarget{example:4.1} \\
\textbf{Example 4.1:} Take $g \equiv 0$. Then the function $h$ takes the form
$$h(t)=1 - \int^{t}_{0} e^{g(\tau)}d\tau = 1 - t$$
and for every $x,y \in \overline{D}(0,1)$:
$$f_g(x,y)= h(x^2+y^2) = 1-(x^2+y^2).$$ \\
\textbf{Example 4.2:} Take $g(t) = -t - \ln(1 - e^{-1})$. Note that $g'(t) \equiv -1 \geq -\frac{1}{t}$ and 
$$\int^{1}_{0} e^{g(t)} dt = \int^{1}_{0} e^{-t - \ln(1 - e^{-1})}dt = \frac{1}{1 - e^{-1}} \int^{1}_{0} e^{-t}dt = \frac{1-e^{-1}}{1 - e^{-1}} = 1.$$
So $g \in \mathcal{G}$. The function $h$ is given by
$$\int^{t}_{0} e^{-t - \ln(1 - e^{-1})}dt = \frac{1}{1 - e^{-1}} \int^{t}_{0} e^{-t}dt = 
\frac{1-e^{-t}}{1 - e^{-1}}$$
and $f_g(x,y)=h(x^2+y^2) = \frac{1-e^{-(x^2+y^2)}}{1 - e^{-1}}$. \\\\
\subsection{Existence of a minimal solution}
In this section we prove that example \hyperlink{example:4.1}{4.1} is special in the following sense: among all solutions $f \in \mathcal{S}$ such that $f(0,0)=1$, the function $f_0=1-(x^2+y^2)$ is minimal.
\begin{Theorem11} \rm{Let $f \in \mathcal{S}$. From theorem \hyperlink{theorem:1.3}{1.3} we know that there exists a unique $g \in \mathcal{S}$ such that $f=f_g$. Denote $h(t) = 1 - \int_0^t e^{g(t)} dt$. Also denote $g_0 \equiv 0$ and $h_0 = 1 - \int_0^t e^{g_0(t)} dt = 1 - \int_0^t 1 dt = 1 - t$. 

First, we want to show that $g(0) \geq g_0(0)=0$. Assume that is not the case. Then there exists some $\epsilon>0$ such that $g(0) \leq - \epsilon$. $g$ is monotone decreasing so $g(t) \leq - \epsilon$ for every $t \in [0,1]$. But then 
$$\int_0^{1} e^{g(t)} dt \leq \int_0^{1} e^{- \epsilon} dt = e^{-\epsilon} < 1 = \int_0^{1} e^{g(t)} dt.$$
This is a contradiction to the assumption. So $g(0) \geq g_0(0)=0$.

The functions $e^{g}$ and $e^{g_0}$ are strictly positive and $e^{g(0)} \geq e^{g_0(0)} = 1$. Since the area under both graphs must equal to $1$ and $e^g$ is monotone decreasing, there exists an intersection point of the graphs. That is, there is $t_0 \in [0,1]$ such that $$\forall t \leq t_0 : e^g(t) \geq e^{g_0(t)} = 1 \: \: \mathrm{and} \: \: \forall t \geq t_0 : e^g(t) \leq e^{g_0(t)} = 1.$$ 
So the graphs look like the following figure:
\begin{figure}[H]
	\centering
	\includegraphics[width=0.6\textwidth]{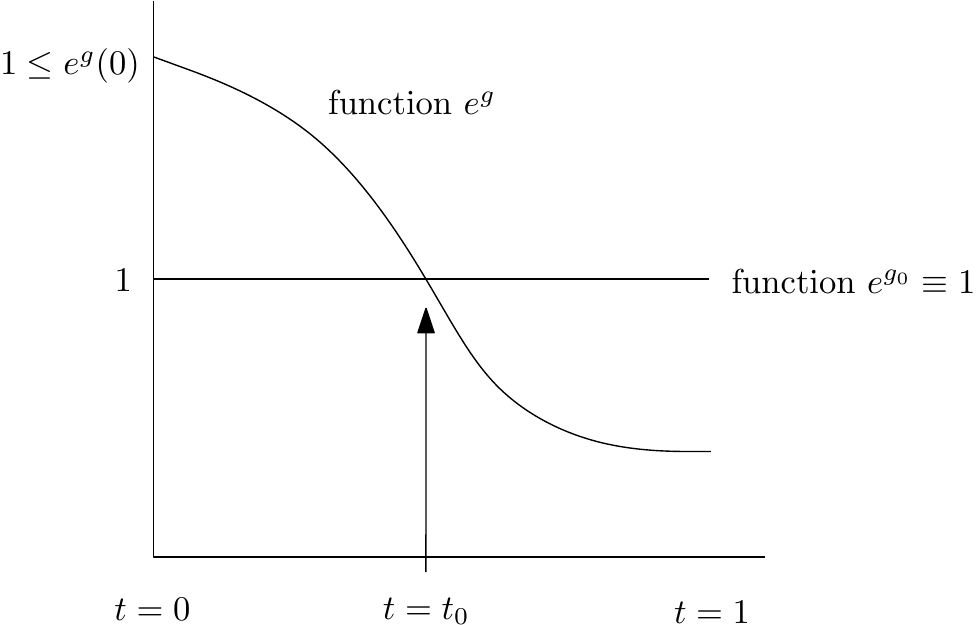}
	\caption{the functions $e^g$ and $e^{g_0}$.}
	\label{fig_2}	
\end{figure}
Note that $\int_0^{1} e^{g(t)} dt = 1$ so the total area under the graph of $g$ must equal to $1$.  
The area under the graph of $e^{g_0} \equiv 1$ obviously equals to $1$ as well. So the areas marked in the figure below must be equal. 
\begin{figure}[H]
	\centering
	\includegraphics[width=0.6\textwidth]{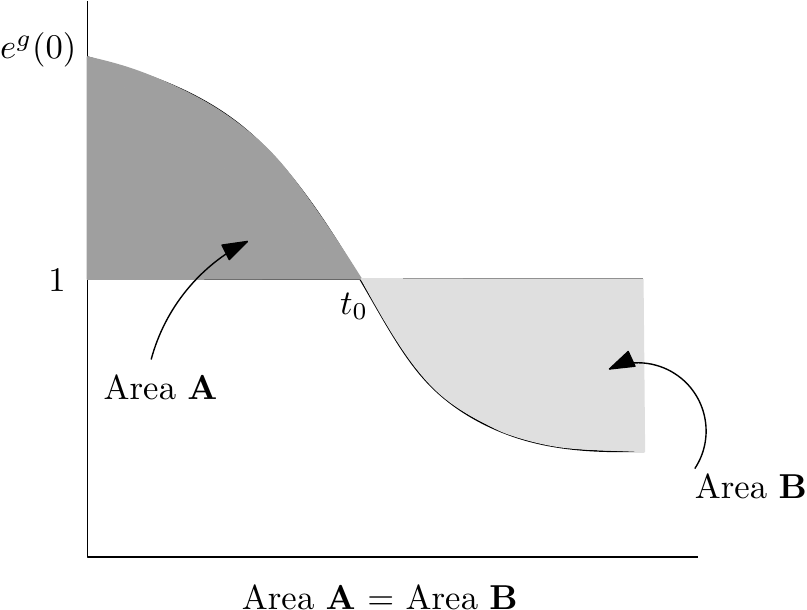}
	\caption{the colored areas \textbf{A} and \textbf{B} are equal.}
	\label{fig_3}	
\end{figure} 

From the figure, it is obvious that for every $t \in [0,1]$:
$$\int_0^t e^{g(t)} dt \geq \int_0^t e^{g_0(t)} dt = \int_0^t 1 dt = t. $$

Thus,

$$h(t) = 1 - \int_0^t e^{g(t)} dt \leq 1 - \int_0^t e^{g_0(t)} dt = h_0(t) = 1 - t. $$

And for every $(x,y) \in \overline{D}(0,1)$:
$$ f(x,y) = h(x^2+y^2) \geq h_0(x^2+y^2) = f_0(x,y) = 1 - (x^2 + y^2).$$
} \end{Theorem11}

\end{document}